\documentclass[11pt,leqno,oneside]{amsart}
\usepackage[width=6.5in,height=8.8in]{geometry}
\usepackage[mathscr]{eucal}
\usepackage{amssymb, amsmath,array, amscd, bbm, comment, charter}
\usepackage{enumerate}
\usepackage{graphicx}
\usepackage{url}
\usepackage[colorlinks,plainpages]{hyperref}

\newtheorem{thm}{Theorem}[section]
\newtheorem{defn}[thm]{Definition}
\newtheorem{lem}[thm]{Lemma}
\newtheorem{cor}[thm]{Corollary}

\theoremstyle{definition}

\begin{document}

\title[Multivolumes of Random Amoebas]{Expected Multivolumes of Random Amoebas}

\author{Turgay Bayraktar}
\thanks{T.\ Bayraktar is partially supported by T\"{U}B\.{I}TAK grant ARDEB-1001/119F184 and Turkish Academy of Sciences, GEBIP}
\address{Faculty of Engineering and Natural Sciences, Sabanc{\i} University, \.{I}stanbul, Turkey}
\email{tbayraktar@sabanciuniv.edu}

\author{Al\.i Ula\c{s} \"Ozg\"ur K\.i\c{s}\.isel}
\thanks{A. U. \"O K\.i\c{s}\.isel is partially supported by CIMPA Research in Pairs Fellowship, T\"{U}B\.{I}TAK 2219 grant and Camille Jordan Institute.} 
\address{Mathematics Department, Middle East Technical University, Ankara, Turkey and Camille Jordan Institute, Lyon, France}
\email{akisisel@metu.edu.tr}
\date{\today}


\begin{abstract}
 We compute the expected multivolume of the amoeba of a random half dimensional complete intersection in $\mathbb{CP}^{2n}$. We also give a relative generalization of our result to the toric case. 
\end{abstract}

\maketitle

\section{Introduction}
Let $n, m$ be non-negative integers such that $n<m$. Suppose that $V$ is an $n$-dimensional algebraic subvariety of $(\mathbb{C}^{\times})^{m}$. The amoeba $\mathcal{A}(V)=Log(V)$ of $V$, as defined in \cite{GKZ}, is the image of $V$ in $\mathbb{R}^m$ under the coordinatewise logarithm map $Log: (\mathbb{C}^{\times})^{m}\rightarrow \mathbb{R}^{m}$:  
\[ Log(z_{1},\ldots,z_{m})=(\log|z_{1}|,\log|z_{2}|,\ldots,\log|z_{m}|). \]
The preimage of a point in $\mathbb{R}^{m}$ under the $Log$ map is an $m$-dimensional (real) torus, and the map $Log$ is a proper map. 

There are several natural measures on amoebas of varieties that have been considered in the literature. In \cite{PR}, Passare and Rullg\r{a}rd consider the case where $V$ is a hypersurface defined by a Laurent polynomial $f$. In this case, there is a naturally defined convex function, the Ronkin function $N_{f}$  \cite{Ron} on $\mathbb{R}^{m}$ which is linear on each connected component of the complement of $\mathcal{A}(V)$. Therefore, the Monge-Amp{\`e}re measure $\mu_{f}$ associated to $N_{f}$ is supported on $\mathcal{A}(V)$. Passare and Rullg\r{a}rd prove (\cite[Theorem 4]{PR}) that the total mass of $\mu_{f}$ is equal to the volume of the Newton polytope of $f$.  If $V$ is a plane curve, then they compare this measure to the Lebesgue measure, and in \cite[Theorem 7]{PR}  they prove that $\mu_{f}$ is greater than or equal to $\pi^{-2}$ times the Lebesgue measure on $\mathcal{A}(V)$. These two results imply that the total Lebesgue area of $\mathcal{A}(V)$ in the case of plane curves is bounded above by $\pi^{2}$ times the area of the Newton polygon of $f$. This bound is sharp, and it is attained for simple Harnack curves (see \cite{Mik1} and \cite{MR}). 

Madani and Nisse \cite{MN} prove that Lebesgue induced volume of a generic $n$-dimensional algebraic subvariety of $(\mathbb{C}^{\times})^{m}$ is finite when $m\geq 2n$. Later in \cite{Mik2}, Mikhalkin considers the case $m=2n$, namely the case where $V$ is a half-dimensional subvariety of $(\mathbb{C}^{\times})^{m}$. In this case, $\dim(\mathcal{A}(V))\leq 2n$, and equality occurs for a generic choice of $V$. If $\dim(\mathcal{A}(V))<2n$, then the restriction of the Lebesgue measure of $\mathbb{R}^{2n}$ to $\mathcal{A}(V)$ is trivial. If $\dim(\mathcal{A}(V))=2n$, then the map $Log|_{V}$ is a finite map, whose degree is bounded above by a toric intersection number as shown in \cite[Theorem 1]{Mik2}. 
Therefore, it makes sense to consider the multivolume measure on $\mathcal{A}(V)$, defined by the formula:  
\[ dMultiVol=Log_{*}(Log^{*}(dVol)),\] 
where $dVol$ denotes the Lebesgue measure on $\mathcal{A}(V)$ induced from $\mathbb{R}^{2n}$. Mikhalkin proves an upper bound for the total multivolume on $\mathcal{A}(V)$ in the case where $V$ is a complete intersection defined by $n$ Laurent polynomials $f_{1},\ldots,f_{n}$. Denote the Newton polytope of $f_{i}$ by $\Delta_{i}$ and let $\alpha(V)$ be the mixed volume of the polytopes $\Delta_{1}, \Delta_{1}, \Delta_{2}, \Delta_{2},\ldots, \Delta_{n}, \Delta_{n}$. In \cite[Lemma 6]{Mik2}, Mikhalkin proves that the total multivolume of $\mathcal{A}(V)$ is bounded above by $\pi^{2n}$ times $\alpha(V)$. As a consequence (\cite[Theorem 1]{Mik2}), the total Lebesgue volume of $\mathcal{A}(V)$ is bounded above by $\pi^{2n}/2$ times $\alpha(V)$. This result coincides with the result of Passare and Rullg\r{a}rd in the case $n=1$. 

In this paper, we study the expected multivolume of a randomly chosen half dimensional complete intersection in $(\mathbb{C}^{\times})^{2n}$. By a random choice of the subvariety we mean the following: Let $W_{d}$ denote the vector space of complex homogeneous degree $d$ polynomials in $2n+1$ variables. There is a $U(2n+1)$-invariant Gaussian probability measure on $W_{d}$ given by the formula 
\[ \frac{1}{\nu_{d}}e^{-||f||^2_{FS}}df \] 
where $||f||_{FS}$ denotes the $L^2$-Fubini-Study norm of $f$, $df$ denotes the Lebesgue measure on $W_{d}$, and $\nu_{d}$ is a normalizing factor. We will describe this measure in more detail in Section 3. This probability distribution is called the Kostlan distribution on $W_{d}$, and statistical properties of solutions of random polynomial systems (especially in the real case) with respect to this distribution have been studied by many authors (for instance, see \cite{EK} and \cite{SS}). Our first result is: 

\begin{thm} \label{main}
Let $d_{1},d_{2},\ldots,d_{n}$ be positive integers. Suppose that $V$ is the subvariety of  $(\mathbb{C}^{\times})^{2n}$ defined by random polynomials $f_{1},\ldots,f_{n}$, where $f_{i}$ is chosen independently  from $W_{d_{i}}$ with respect to the Kostlan distribution. Then, the expected total multivolume of $\mathcal{A}(V)$ is given by the formula: 
\[  \mathbb{E}(MultiVol(\mathcal{A}(V)))= \pi^{2n} d_{1}d_{2}\ldots d_{n}. \]
\end{thm} 
In the case of plane curves, we note that this expected value is equal to $\pi^2 d$, which is significantly smaller than the maximum possible value of $\pi^2 d^2$ as $d\rightarrow +\infty$. 

The organization of the paper is as follows: In section 2, we review rolled coamoebas and an essential property of the amoeba and rolled coamoeba maps. In section 3, we carefully define random complete intersections with respect to Kostlan distribution. Section 4 contains the proof of Theorem \ref{main}. In section 5, we discuss a toric generalization of our theorem: More precisely, by using a result of Malajovich \cite{Mal}, we prove a relative generalization of Theorem \ref{main} to the toric case.  

\section*{Acknowledgement}
The authors thank Michele Ancona, Lie Fu, Ilia Itenberg, Antonio Lerario, Grigory Mikhalkin and Jean-Yves Welschinger for helpful discussions. The authors thank Antonio Lerario for bringing the reference \cite{Mal} to their attention.

\section{Amoebas and Rolled Coamoebas} 
In this section, we will recall some definitions and results from \cite{Mik2} that will be used in the proofs of our theorems. Let $V$ be a complex $n$-dimensional algebraic subvariety of $(\mathbb{C}^{\times})^{2n}$. The rolled coamoeba of $V$ is the image of $V$ under the ``coordinatewise argument modulo $\pi$ map'' $Arg_{\pi}: (\mathbb{C}^{\times})^{2n} \rightarrow \mathbb{T}$: 
\[ Arg_{\pi}(z_{1},\ldots,z_{2n})=(\arg(z_{1}),\ldots,\arg(z_{2n})) \] 
where $\mathbb{T}=(\mathbb{R}/\pi \mathbb{Z})^{2n}\cong (S^{1}_{\pi})^{2n}$. 
Notice that $Log+iArg_{\pi}$ locally agrees with the coordinatewise holomorphic logarithm function on $(\mathbb{C}^{\times})^{2n}$. 

As opposed to the $Log$ map, $Arg_{\pi}$ is not in general a proper map when restricted to $V$. For example, if $n=1$ and $V$ is the line given by the equation $x+y=1$ in $(\mathbb{C}^{\times})^{2}$ , then the preimage $Arg_{\pi}^{-1}(0,0)$ coincides with the set of real solutions of this equation inside $(\mathbb{C}^{\times})^{2}$, which is homeomorphic to $\mathbb{R}-\{0,1\}$. If $\Theta\in \mathbb{T}$ is a regular value of $Arg_{\pi}$ in its image, then by the inverse function theorem $Arg_{\pi}$ is locally a homeomorphism in a neighborhood in $V$ of any point in $Arg_{\pi}^{-1}(\Theta)$ to a neighborhood of $\Theta$ in $\mathbb{T}$. By Sard's theorem, the set of critical values of $Arg_{\pi}$ has Lebesgue measure zero in $\mathbb{T}$. Therefore, there are two alternatives: Either the image $Arg_{\pi}(V)$ has real dimension less than $2n$, or this dimension is equal to $2n$ and the image contains a dense set of regular values. Both of these alternatives can occur, while for a generic choice of $V$ the second one takes place. As an example for the first alternative, one can consider the line $x=y$ in $(\mathbb{C}^{\times})^{2}$. Then the  image $Arg_{\pi}(V)$ will be the diagonal subset in $\mathbb{T}$. 

The observation that $Log$ and $Arg_{\pi}$ are the real and imaginary parts of a holomorphic function on the half-dimensional algebraic subvariety $V$ has an important consequence which we will describe now. The following Lemma and its Corollary are paraphrased versions of Lemma 4 and Corollary 5 in \cite{Mik2}, therefore we state them without proof. 

\begin{lem} 
Suppose that $p\in V$ and $U$ is an open neighborhood of $p$ in $V$ such that both of the maps $Log$ and $Arg_{\pi}$ are homeomorphisms from $U$ to their respective images. Then the map 
\[  Log\circ Arg_{\pi}^{-1} : Arg_{\pi}(U)\rightarrow Log(U) \] 
is volume preserving, with respect to the Lebesgue measures on $\mathbb{T}$ and $\mathbb{R}^{2n}$. 
\end{lem}

\begin{cor} \label{preserve}
Let $d\Theta$ denote the Lebesgue measure on $\mathbb{T}$ restricted to the coamoeba of $V$. Then 
\[ Log_{*}(Arg_{\pi}^{*}(d\Theta))=dMultiVol. \] 
\end{cor}

\section{Random Complex Kostlan Polynomials} 
Our next goal is to make a careful definition of random  complete intersections in $(\mathbb{C}^{\times})^{m}$ with respect to the Kostlan distribution. At this stage, we will have no restictions on the Newton polytope except the degrees of the defining polynomials, hence the process can be viewed as choosing random  complete intersections in $\mathbb{CP}^{m}$, assuming that  $(\mathbb{C}^{\times})^{m}$ is identified with the complement of the coordinate hyperplanes in the affine chart $X_{0}\neq 0$ in $\mathbb{CP}^{m}$. 
In the last section, we will discuss the toric generalizations, namely the case of restricted Newton polytopes. An exposition for the following discussion can be found in \cite[sections 2.1-2.4]{Wel}. 

Let $W_{d}=\mathbb{C}_{d}^{hom}[X_{0},\ldots,X_{m}]$ denote the space of homogeneous polynomials of degree $d$ in $m+1$ variables. Then $W_{d}$ can be identified with the space of sections $H^{0}(\mathbb{CP}^{m},\mathcal{O}(d))$, where $\mathcal{O}(d)=\mathcal{O}(1)^{\otimes d}$ and $\mathcal{O}(1)$ denotes Serre's twisting sheaf. The complex projective space $\mathbb{CP}^{m}$ can be viewed as a homogeneous space 
\[  \mathbb{CP}^{m} \cong U(m+1)/U(1)\times U(m) \] 
and it inherits a unitary invariant metric $h$ from $U(m+1)$, which is the well-known Fubini-Study metric. We assume that the Fubini-Study metric is normalized so that its total volume on $\mathbb{CP}^{m}$ is equal to $1$. The metric $h$ naturally induces pointwise Hermitian products $h_{d}$, and in turn induces  $L^{2}$-Hermitian products 
\[  \langle f, g\rangle = \int_{\mathbb{CP}^{m}} h_{d}(f,g) dVol_{FS}  \] 
on each $W_{d}$. Here, $dVol_{FS}$ is the Fubini-Study volume form. Monomials $X_{0}^{\alpha_{0}}X_{1}^{\alpha_{1}}\ldots X_{m}^{\alpha_{m}}$ form an orthogonal basis for the $L^{2}$-Hermitian space $W_{d}$. Their squared norms are given by the inverse multinomial coefficients
\[  ||X_{0}^{\alpha_{0}}X_{1}^{\alpha_{1}}\ldots X_{m}^{\alpha_{m}}||^{2}_{FS}= \left(\frac{d!}{\alpha_{0}! \alpha_{1}!\ldots,\alpha_{m}!}\right)^{-1}. \]

The complex vector space $W_{d}$ can now be equipped by the following Gaussian probability distribution, called the Kostlan distribution 
\[  d\mu_{\mathbb{C}}(f)=\frac{1}{\pi^{N_{d}}} e^{-||f||^{2}_{FS}}df. \] 
Here, $N_{d}$ denotes the (complex) dimension of the vector space $W_{d}$ and $df$ denotes the Lebesgue measure on $W_{d}$. 

Each complex polynomial can be decomposed uniquely into a sum of a polynomial with real coefficients and a polynomial with purely imaginary coefficients. Hence, we obtain a direct sum decomposition $W_{d}=Re(W_{d})\oplus i Im(W_{d})$. The real vector spaces $Re(W_{d})$ and $Im(W_{d})$ inherit probability distributions $d\mu_{\mathbb{R}}(f)$ and $d\mu_{i\mathbb{R}}(f)$ by a similar formula, with only the normalizing constant modified: 
\[ d\mu_{\mathbb{R}}(f)=d\mu_{i\mathbb{R}}(f)=\frac{1}{\sqrt{\pi}^{N_{d}}} e^{-||f||^{2}_{FS}}df. \]
The measure $d\mu_{\mathbb{C}}(f)$ is in fact the product measure of $d\mu_{\mathbb{R}}(f)$ and $d\mu_{i\mathbb{R}}(f)$. 

\begin{defn} 
Let $k$ and $d_{1},d_{2},\ldots,d_{k}$ be positive integers. A random complete intersection of degree sequence $(d_{1},d_{2},\ldots,d_{k})$ in $\mathbb{CP}^{m}$ is the common zero locus of $f_{1},\ldots,f_{k}$ where $(f_{1},\ldots,f_{k})$ is chosen independently at random from the sample space $W_{d_{1}}\oplus W_{d_{2}} \oplus \ldots \oplus W_{d_{k}}$ equipped with the product probability measure of the Kostlan measures described above. 
\end{defn} 
We remark that a random complete intersection will almost surely have codimension $k$ and hence it will be an actual complete intersection in the usual sense. 

\section{Proof of Theorem \ref{main}} 
We first recall a theorem of Shub and Smale on the expected number of common real solutions of a full random real polynomial system of equations where the polynomials are chosen with respect to Kostlan distribution. In accordance with the notation of the previous section, let $k=m$ and let $f_{1},\ldots,f_{k}$ be real polynomials of degrees $d_{1},\ldots,d_{k}$ respectively such that $f_{i}$ is randomly chosen from $Re(W_{d_{i}})$ with its measure described as in the previous section. The theorem of Shub and Smale (\cite[Theorem A]{SS}) is as follows: 

\begin{thm}  \label{ShubSmale}
The expected number of common real zeros of $f_{1},f_{2},\ldots,f_{k}$ in $\mathbb{RP}^{k}$ is $(d_{1}d_{2}\ldots d_{k})^{1/2}.$ 
\end{thm} 
We remark that the common real zeros will almost surely lie in $(\mathbb{R}^{\times})^{k}$. Next, we prove our main result: 

\noindent \textit{Proof of Theorem \ref{main}:}  Let $m=2n$, $k=n$ in the notation of section 3. Note that, by Corollary \ref{preserve}
\begin{eqnarray*}
 MultiVol(\mathcal{A}(V))&=&\int_{\mathcal{A}(V)}dMultiVol \\
&=& \int_{\mathcal{A}(V)}Log_{*}(Arg_{\pi}^{*}(d\Theta)) \\ 
&=& \int_{V} Arg_{\pi}^{*}(d\Theta) \\
&=& \int_{\mathbb{T}} \#(Arg_{\pi}^{-1}(\Theta))d\Theta 
\end{eqnarray*} 
where $ \#(Arg_{\pi}^{-1}(\Theta))$ denotes the number of preimages of $\Theta$ in $V$. On the other hand, by using Fubini's theorem, we get 
\[  \mathbb{E}(MultiVol(\mathcal{A}(V))=\int_{\mathbb{T}} \mathbb{E}( \#(Arg_{\pi}^{-1}(\Theta))) d\Theta. \] 
Since $\int_{\mathbb{T}} d\Theta = \pi^{2n}$, the proof of the theorem will be complete if we prove that  $$\mathbb{E}( \#(Arg_{\pi}^{-1}(\Theta)))=d_{1}d_{2}\ldots d_{n}$$ for every $\Theta\in \mathbb{T}$.  

Let us first prove that for $\Theta=(0,0,\ldots,0)=\mathbf{0}$, we have $\mathbb{E}( \#(Arg_{\pi}^{-1}(\mathbf{0})))=d_{1}d_{2}\ldots d_{n}$. Each of the complex coordinates $(z_{1},z_{2},\ldots,z_{2n})$ of $(\mathbb{C}^{\times})^{2n}$ can be written uniquely as $z_{j}=r_{j}e^{i\theta_{j}}$ for $r_{j}\in \mathbb{R}\setminus \{0\}$ and $\theta_{j}\in [0,\pi)$. Given that $\Theta=\mathbf{0}$, for each $j$ we have $\theta_{j}=0$. Therefore, $\#(Arg_{\pi}^{-1}(\mathbf{0}))$ is equal to the number of common real zeros of the complex polynomials $f_{i}(r_{1},r_{2},\ldots,r_{2n})$ for $i=1,2,\ldots,n$. Each $f_{i}$ has a real and imaginary part which are both real polynomials of degree $d_{i}$ and which must both vanish whenever $f_{i}$ does. Furthermore, the complex Kostlan distribution on $W_{d_{i}}$ induces the real Kostlan distribution on $Re(W_{d_{i}})$ and $Im(W_{d_{i}})$. Therefore, applying Theorem \ref{ShubSmale} to the polynomials $Re(f_{1}), Im(f_{1}),\ldots,Re(f_{n}),Im(f_{n})$, we obtain 
\[ \mathbb{E}( \#(Arg_{\pi}^{-1}(\mathbf{0})))=(d_{1}d_{1}d_{2}d_{2}\ldots d_{n}d_{n})^{1/2}=d_{1}d_{2}\ldots d_{n},\]
proving the claim in this case. 

For the general case, notice that the action of the unitary group $U(2n+1)$ on $\mathbb{CP}^{2n}$ by the formula $g\cdot \mathbf{X}= g\mathbf{X}$, where $\mathbf{X}=[X_{0}, X_{1},\ldots, X_{2n}]^{t}$,  has an induced action on $W_{d}$ by $P\rightarrow P\circ g^{-1}$. This action gives us an isometry of $(W_{d}, h_{d})$ for every $g\in U(2n+1)$. We note that this action does not leave $Re(W_{d})$ and $Im(W_{d})$ invariant, but we still have the orthogonal direct sum decomposition $W_{d}=g\cdot Re(W_{d})\oplus g\cdot Im(W_{d})$ and $d\mu_{\mathbb{C}(f)}$ is the product measure of the Gaussian measures on these subspaces. Now, given any $\Theta=(\theta_{1},\ldots,\theta_{2n})\in \mathbb{T}$, there exists a diagonal matrix $g=diag(1,e^{-i\theta_{1}},\ldots,e^{-i\theta_{2n}})$ in $U(2n+1)$ such that for any $\mathbf{X}$ with $Arg_{\pi}(\mathbf{X})=\Theta$, we have $Arg_{\pi}(g\cdot \mathbf{X})=\mathbf{0}$. The corresponding isometry on $W_{d}$ will leave the expected number of preimages of $Arg_{\pi}^{-1}$ invariant. Hence 
$\mathbb{E}( \#(Arg_{\pi}^{-1}(\Theta)))=\mathbb{E}( \#(Arg_{\pi}^{-1}(\mathbf{0})))=d_{1}d_{2}\ldots d_{n}.$  $ \Box$ 

\begin{cor}
Under the assumptions of Theorem \ref{main}, the expected Lebesgue volume of $\mathcal{A}(V)$ is bounded above by $\pi^{2n}d_{1}\ldots d_{n} /2$. 
\end{cor} 

\noindent \textit{Proof}  The $Log$ map has degree at least $2$ on a generic point of $\mathcal{A}(V)$ (see, for instance, \cite{Mik2}, Corollary 7). The claim follows. $\Box$

\section{Toric Generalization} 
Computation of simultaneous complex zeros of deterministic as well as Gaussian systems of sparse polynomials has been studied by various authors (see eg. \cite{HSt,Rojas,MaR,SZ}) by using methods of algebraic and toric geometry. More recently, the first author proved that limiting distribution of (complex zeros) of random Laurent polynomials whose support are contained in dilates of a fixed integral polytope is determined by associated pluri-potential theory \cite{Bay}.

A natural question is how one can generalize Theorem \ref{main} to the case where each of the defining polynomials $f_{i}$ is sparse, namely, chosen at random such that it has Newton polytope $\Delta_{i}$. The analogous generalization of the theorem of Shub and Smale in this case for real zeros is not completely clear. For some developments in this direction, see \cite{MaR} and \cite{BETC}. 

However, in \cite{Mal}, Theorem 1.3 and its corollary, Malajovich obtained a relative toric generalization of the Shub-Smale theorem in the sense that each $f_{i}$ is chosen to have a scaled polytope $d_{i}\Delta_{i}$, and the probability measure on the toric variety for the scaled polytope $d_{i}\Delta_{i}$ is obtained from the one on $\Delta_{i}$ in a natural way. The resulting probability measures also turn out to be invariant under the diagonal subgroup of the unitary group, hence this allows us to generalize Theorem \ref{main} to this set-up. We will describe this in detail now. The facts stated in the paragraphs below on toric varieties can be found in the first two chapters of the book \cite{CLS}.   

Let $M=\mathbb{Z}^{n}$ and $M_{\mathbb{R}}=M\otimes \mathbb{R}$. Let $\Delta\subset M_{\mathbb{R}}$ be a full-dimensional normal lattice polytope (see \cite{CLS}, definition 2.2.9). Let $A=\Delta\cap M$.  Every $m=(a_{1},\ldots,a_{m})\in A$ defines a character $\chi^{m}: (\mathbb{C}^{\times})^{n}\rightarrow \mathbb{C}^{\times}$ by
\[ \chi^{m}(t_{1},\ldots,t_{m})=t_{1}^{a_{1}}\ldots t_{m}^{a_{m}}. \] 
In turn, we obtain a map $\Phi_{\Delta}: (\mathbb{C}^{\times})^{n} \rightarrow \mathbb{P}(\mathbb{C}^{A})$ by 
\[ \Phi_{\Delta}=(\chi^{m_{1}},\ldots,\chi^{m_{s}}) \] 
where $A=\{m_{1},\ldots,m_{s}\}$. The Zariski closure of $ \Phi_{\Delta}((\mathbb{C}^{\times})^{n})$ inside $\mathbb{P}(\mathbb{C}^{A})$ is called the projective toric variety associated to $\Delta$, and is denoted by $X_{\Delta}$. This construction not only associates the abstract variety $X_{\Delta}$ to $\Delta$, but also a preferred embedding into projective space, i.e. a choice of line bundle. For example, if $\Delta$ is a simplex of side-length $d$ (a $d$-times scaled version of the unit simplex), then we obtain $X_{\Delta}=\mathbb{P}^{n}$ with its $d$-uple (Veronese) embedding into projective space. 

Let us denote the embedding of $X_{\Delta}$ into $\mathbb{P}(\mathbb{C}^{A})$  discussed in the above paragraph by $i:X_{\Delta}\rightarrow \mathbb{P}(\mathbb{C}^{A})$. Then, the pull-back of the twisting sheaf $\mathcal{O}(1)$ on the target space $\mathbb{P}(\mathbb{C}^{A})$ to $X_{\Delta}$ induces the complex spaces of sections 
\[ W_{\Delta}=H^{0}(X_{\Delta}, i^{*}(\mathcal{O}(1))), \] 
which generalize the spaces $W_{d}$ from section 3. The set of characters $\chi^{m}$ above can be identified with a basis of $W_{\Delta}$. Let $h_{\Delta}$ be any Hermitian metric on $W_{\Delta}$ such that the characters $\chi^{m}$ are orthogonal. Let $d\Delta$ denote the Minkowski sum of $\Delta$ with itself $d$ times. Notice that the sheaf giving us $W_{d\Delta}$ is the $d$-th tensor product of the one for $W_{\Delta}$. Then, $h_{\Delta}$ naturally induces a Hermitian metric $h_{d\Delta}$ on $W_{d\Delta}$ for every positive integer $d$; the metric $h_{d\Delta}$ agrees with the one obtained by Aronszajn multiplication formula of reproducing kernel spaces (see \cite{Mal}). Furthermore, with respect to $h_{d\Delta}$, the characters associated to the polytope $d\Delta$ will remain orthogonal. On each $W_{\Delta}$, the metric $h_{\Delta}$ defines an $L^2$-Hermitian product by the formula
\[ \langle f,g \rangle=\int_{X_{\Delta}} h_{\Delta}(f,g) dVol_{\Delta}. \] 
In turn, we obtain a norm, and a Gaussian probability measure on $W_{\Delta}$: 
\[ d\mu_{\Delta} = \frac{1}{\nu_{d}} e^{-||f||^{2}}df.  \] 
where $\nu_{d}$ is a normalizing constant. 

Suppose now that $(f_{1},f_{2},\ldots,f_{n})$ is a randomly chosen $n$-tuple of complex polynomials from the probability space  $W_{\Delta_{1}}\oplus \ldots \oplus W_{\Delta_{n}}$ with its product probability measure. Denote the expected multivolume of the amoeba of the complete intersection determined by these polynomials in $(\mathbb{C}^{\times})^{2n}$ by $\mathbb{E}_{\Delta_{1},\ldots,\Delta_{n}}$. 

\begin{thm} 
Let $d_{1},d_{2},\ldots,d_{n}$ be positive integers. Then, 
\[  \mathbb{E}_{d_{1}\Delta_{1},\ldots,d_{n}\Delta_{n}}=(d_{1}d_{2}\ldots d_{n}) \mathbb{E}_{\Delta_{1},\ldots,\Delta_{n}}. \] 
\end{thm} 

\noindent \textit{Proof} The proof goes along the lines of the proof of Theorem \ref{main} with the following addenda: Instead of the theorem of Shub-Smale, we use the theorem of Malajovich (\cite{Mal}, Theorem 1.3 and Corollary 1.4). Also, we do not necessarily have full unitary invariance of the metric. However, since the orthogonality of the characters were assumed, the invariance under the diagonal subgroup of the unitary group is retained, and this suffices for the proof. $\Box$

\end{document}